\documentclass[10pt]{article}
\usepackage{amssymb,latexsym,amsmath,epsfig,amsthm}

 \makeatletter

\renewcommand\section{\@startsection {section}{1}{\z@}
{-30pt \@plus -1ex \@minus -.2ex}
{2.3ex \@plus.2ex}
{\normalfont\normalsize\bfseries}}

\renewcommand\subsection{\@startsection{subsection}{2}{\z@}
{-3.25ex\@plus -1ex \@minus -.2ex}
{1.5ex \@plus .2ex}
{\normalfont\normalsize\bfseries}}

\renewcommand{\@seccntformat}[1]{\csname the#1\endcsname. }

\makeatother

\begin{document}

\begin{center}
\uppercase{\bf Multiplicative bases and an Erd\H{o}s problem}
\vskip 20pt
{\bf P\'eter P\'al Pach \footnote{{Department of Computer Science and Information Theory, Budapest University of Technology and Economics, 1117 Budapest, Magyar tud\'osok k\"or\'utja 2., Hungary}\\
{\tt ppp@cs.bme.hu}. This author was supported by the Hungarian Scientific Research Funds (Grant Nr. OTKA PD115978 and OTKA K108947) and the J\'anos Bolyai Research Scholarship of the Hungarian Academy of Sciences.}, Csaba S\'andor \footnote{{Institute of Mathematics, Budapest University of
Technology and Economics, H-1529 B.O. Box, Hungary}\\ {\tt csandor@math.bme.hu}.
This author was supported by the OTKA Grant No. K109789. This paper was
supported by the J\'anos Bolyai Research Scholarship of the Hungarian
Academy of Sciences. }}\\

\bigskip

\end{center}

\noindent

\pagestyle{myheadings}
\thispagestyle{empty}
\baselineskip=12.875pt

\newtheorem*{rado}{Rado's Theorem}
\newtheorem *{sarkozy2}{Theorem B (S\'ark\"ozy, \cite{Sarkozy2})}
\newtheorem *{sarkozy1}{Theorem A (S\'ark\"ozy, \cite{Sarkozy1})}
\newtheorem{theorem}{Theorem}
\newtheorem{lemma}[theorem]{Lemma}
\newtheorem{problem}{Problem}
\newtheorem{corollary}[theorem]{Corollary}
\newtheorem{example}{Example}[section]
\newtheorem{proposition}[theorem]{Proposition}

\begin{abstract}
In this paper we investigate how small the density of a multiplicative basis of order $h$ can be in $\{1,2,\dots,n\}$ and in $\mathbb{Z}^+$. Furthermore, a related problem of Erd\H os is also studied: How dense can a set of integers be, if none of them divides the product of $h$ others?
\end{abstract}


\section{Introduction}
Throughout the paper we are going to use the notions $[n]=\{1,2,\dots,n\}$ and $A(n)=A\cap [n]$ for $n\in\mathbb{N}$ and $A\subseteq \mathbb{Z}^+$.
Let $h\ge 2$ and $S\subset \mathbb{Z}^+$. We say that the set $B\subseteq \mathbb{Z}^+$ forms a multiplicative basis of $S$, if every element of $s\in S$ can be written as the product of $h$ members of $B$. The set of these  multiplicative bases will be denoted by $MB_h(S)$.  While the additive basis is a popular topic in additive number theory, much less attention was devoted to the multiplicative basis. It is easy to see that every multiplicative basis $B\in MB_h([n])$ contains the prime numbers up to $n$. Let $G_h(n)$ denote the smallest possible size of a basis in $MB_h([n])$. Chan \cite{chan} proved that there exists  some $c_1>0$ such that for every $h\ge 2$ we have $|G_h(n)|\leq \pi (n)+c_1(h+1)^2\frac{n^{\frac{2}{h+1}}}{\log ^2 n}$ (in fact, he did not use the terminology multiplicative basis). In the first theorem we determine the order of magnitude of $G_h(n)-\pi(n)$ in the sense that $h$ is not fixed, the only restriction is that $n$ has to be large enough compared to $h$.

\begin{theorem}\label{MBveges}
Let $h,n\in\mathbb{Z}^+$ such that $h\leq \sqrt{\frac{\log n}{12\log\log n}}$. 
Then for the smallest possible size of a multiplicative basis of order $h$ for $[n]$ we have
$$\pi (n)+0.5h\frac{n^{2/(h+1)}}{\log ^2n}\le G_h(n)\le \pi (n)+150.4h\frac{n^{2/(h+1)}}{\log ^2n}.$$
\end{theorem}

Raikov \cite{rai} proved in 1938 that a $B\in MB_h(\mathbb{Z}^+)$ must be dense sometimes. 

\begin{theorem}{(Raikov, 1938)}
\label{raikov}
Let $B\in MB_h(\mathbb{Z}^+)$. Then $$\displaystyle \limsup _{n\to \infty} \frac{|B(n)|}{n/\log ^{\frac{h-1}{h}} n }\ge \Gamma \left( \frac{1}{h}\right) ^{-1}.$$
On the other hand,  for every $h\ge 2$ there exists a $B_h\in MB_h(\mathbb{Z}^+)$ such that $\displaystyle \limsup _{n\to \infty} \frac{|B(n)|}{n/\log ^{\frac{h-1}{h}} n }< \infty $.
\end{theorem}

In the lower bound the quantity $\Gamma \left( \frac{1}{h} \right)$ is asymptotically $h$. Our next theorem determines the previous limit superior for multiplicative bases of order $h$ up to a constant factor (not depending on $h$).

\begin{theorem}
\label{rai2}
Let $B\in MB_h(\mathbb{Z}^+)$. Then
$$\displaystyle \limsup _{n\to \infty} \frac{|B(n)|}{n/\log ^{\frac{h-1}{h}} n }\ge \frac{\sqrt{6}}{e\pi }.$$
 On the other hand, there exists some  $C>0$ such that for every $h\ge 2$ one can find a $B_h\in MB_h(\mathbb{Z}^+)$ such that $\displaystyle \limsup _{n\to \infty} \frac{|B(n)|}{n/\log ^{\frac{h-1}{h}} n }=C $.
\end{theorem}

On the other hand, a set $B\in MB_h(\mathbb{Z}^+)$ may be thin as our following theorem shows:

\begin{theorem}
\label{MBliminf}
Let $1<h\in\mathbb{Z}^+$.
If $B\in MB_h(\mathbb{Z}^+)$, then $\liminf\limits_{n\to \infty}\frac{|B(n)|}{{\frac{n}{\log n}}}>1$.
On the other hand, for every $\varepsilon>0$ there exists a $B\in MB_h(\mathbb{Z}^+)$ such that
$$\liminf_{n\to\infty} \frac{|B(n)|}{{\frac{n}{\log n}}}<1+\varepsilon.$$
\end{theorem}

The logarithmic density of a set $B\subset \mathbb{Z}^+$ is defined as the limit  $\displaystyle \lim_{n\to \infty }\frac{\sum\limits_{b\in B(n)}\frac{1}{b}}{\log n}$ (if it exists). Our following theorem determines the possible lower densities of the quantity $\sum\limits_{b\in B(n)}\frac{1}{b}$ for a $B\in MB_h(\mathbb{Z}^+)$.

\begin{theorem}
\label{MBlog}
Let $h\ge 2$ and $B\in MB_h(\mathbb{Z}^+)$. Then
$$ \liminf _{n\to \infty}\frac{\sum\limits_{b\in B(n)}\frac{1}{b}}{h\sqrt[h]{\log n} }\geq \frac{\sqrt{6}}{e\pi }.$$
On the other hand, there exists a constant $C$ such that for every $h\geq 2$ there exists a $B\in MB_h(\mathbb{Z}^+)$ such that
$$\limsup _{n\to \infty}\frac{\sum\limits_{b\in B( n)}\frac{1}{b}}{h\sqrt[h]{\log n} }<C.$$
\end{theorem}

If one looks at the paper \cite{ep} of Erd\H os, it seems that he deals with a quite different problem. However, by a closer look it turns out that his problem is closely related to the multiplicative bases.
We say that $A\subset S$ possesses property $P_h$, if there are no distinct elements $a,a_1,\dots ,a_h\in A$ with $a$ dividing the product $a_1\dots a_h$. Denote the set of these $A$'s by $P_h(S)$. Let $\displaystyle F_h(n)=\max _{A\in P_h([n])}|A|$. Clearly the set of prime numbers satisfies property $P_h$, therefore $F_h(n)\ge \pi (n)$. The case $h=2$, that is, such sets of integers where none of the elements divides the product of two others, was settled by Erd\H os \cite{ep}. Chan, Gy\H ori and S\' ark\"ozy \cite{cgs} studied the case $h=3$. Furthermore, recently Chan \cite{chan} determined the order of magnitude of $F_h(n)-\pi (n)$ for every fixed $h$.

\begin{theorem} (Chan, 2011)
There exist absolute constants $c_2,c_3>0$ such that, for any positive integers $n>e^{48}$ and $2\le h\le \frac{1}{6}\sqrt{\frac{\log n}{\log \log n}}$,
$$\pi (n)+\frac{c_2}{(h+1)^2}\frac{n^{2/(h+1)}}{\log ^2n}\le F_h(n)\le \pi (n)+c_3(h+1)^2\frac{n^{2/(h+ 1)}}{\log ^2n}.$$
\end{theorem}

Our next theorem provides a better estimation for $F_h(n)$. Here, the "error term" in the lower and upper bounds differ only by a constant factor not depending on $h$.

\begin{theorem}
\label{fhveges}
Let $h,n\in\mathbb{Z}^+$ such that $h\leq \sqrt{\frac{\log n}{12\log\log n}}$. Then
$$\pi(n)+0.2\frac{n^{2/(h+1)}}{\log ^2n}\leq F_h(n)\leq \pi(n)+379.2\frac{n^{2/(h+1)}}{\log ^2n}.$$
\end{theorem}

Our following two results show us that a sequence $A\in P_h(\mathbb{Z}^+)$ must be thin sometimes, but it may be as dense as allowed by the obtained upper bound in the finite case.

\begin{theorem}
\label{Erdinf}
Let $2\leq h\in\mathbb{Z}$ and $A\in P_h(\mathbb{Z}^+)$. Then for every $\varepsilon >0$
 $$\liminf _{n\to \infty} \frac{A(n)-\pi (n)}{n^{\varepsilon}}<\infty .$$
On the other hand, there exists a constant $c>0$ such that for every $h\ge 2$ a set $A\in P_h(\mathbb{Z}^+)$ can be constructed in such a way that  $|A(n)|\geq \pi (n)+\exp\left\{{(\log n)^{1-\frac{c\sqrt{\log h}}{\sqrt{\log \log n}}}}\right\}$ holds for every $n$.
\end{theorem}

\begin{proposition}
For every $h\ge 2$ there exists an $A_h\in P_h(\mathbb{Z}^+)$ such that $$\limsup _{n\to \infty}\frac{|A_h(n)|-\pi (n)}{\frac{n^{2/(h+1)}}{\log ^2n}}>0.$$
\end{proposition}
The proof of this proposition is going to be omitted because the construction can be easily built up by repeating  the construction of the finite case for bigger and bigger blocks.

Finally, let us mention that the logarithmic density of a set in $P_h(\mathbb{Z}^+)$ can be easily treated because the prime numbers imply that for every $A\in P_h(\mathbb{Z}^+)$ we have $\displaystyle \sum _{a\in A(n)}\frac{1}{a} > \log \log n-c$. On the other hand, by Theorem~\ref{fhveges} we have for every $A\in P_h{(\mathbb{Z}^+)}$
\begin{multline*}
\sum _{a\in A(n)}\frac{1}{a}=\sum _{k\le n}\frac{|A(k)|-|A(k-1)|}{k}=\sum _{k\le n-1}\frac{|A(k)|}{k(k+1)}+\frac{|A(n)|}{n}\le \\
\sum _{k\leq n-1}\frac{\pi (k)+C_hk^{2/3}}{k(k+1)}+1<\log \log n +c_h.
\end{multline*}

The main part of the paper is organized as follows. In Section 2. we prove Theorems \ref{MBveges} and \ref{fhveges} about the finite case and Section 3. contains the proofs of the results about the infinite case.

\section{Finite case}

At first it is going to be considered how small a multiplicative basis of order $h$ for $[n]$ can be.
During the calculations the following well-known estimates  \cite{Rosser} are going to be used:
\begin{lemma}\label{primbecs}
For every $x\geq17$ we have $\frac{x}{\log x}< \pi(x)$.
For every $x>1$ we have $\pi(x)\leq 1.26\frac{x}{\log x}$.
\end{lemma}

\noindent
Now the proof of Theorem~\ref{MBveges} is going to be presented.

\smallskip
\noindent
{\bf Proof of Theorem \ref{MBveges}}.
Let $n^{\frac{1}{h+1}}(\log n)^{-1}=s$.
We start by proving the first statement. Let us assume that $B$ is a multiplicative basis of order $h$ for $[n]$. Clearly, all the prime numbers not greater than $n$ (and 1) have to be in $B$. Our aim is to show that there are at least  $hs^2/2$ elements in $B$ that are the product of at least two primes. Let $V$ denote the set of primes not greater than $n^{1/(h+1)}$:
$V=\{p\ |\ p\leq n^{1/(h+1)}\text{ and }p\text{ is a prime}\}$. According to Lemma~\ref{primbecs}, the size of $V$ is at least $(h+1)s$. If $\{p_1,p_2,\dots,p_{h+1}\}$ is an $(h+1)$-element subset of $V$, then $a=p_1p_2\dots p_{h+1}\leq n$, so $a\in B^h$ implies that there exists a subset $H$ of $\{p_1,p_2,\dots,p_{h+1}\}$ containing at least 2 elements such that $\prod\limits_{p_i\in H} p_i \in B$.
Let $G$ be the hypergraph with vertex set $V$ and edge set $\mathcal H$, where $\mathcal H$ contains those at least 2-element subsets $H$ of $V$ for which $\prod\limits_{p_i\in H} p_i \in B$.
We have already seen that each $(h+1)$-element subset of $V$ contains at least one hyperedge of $\mathcal H$. As $|B|\geq \pi(n)+|\mathcal H|$, our aim is to give a lower bound for $|\mathcal H|$. If  each set in $\mathcal H$ is replaced by one of its 2-element subsets -- the new set of subsets is denoted by $\mathcal H'$ --, then it still remains true that each $(h+1)$-element subset of $V$ contains an element of $\mathcal H'$. Moreover, $|\mathcal H|\geq |\mathcal H'|$. Let $G'$ be the graph with vertex set $V$ and edge set $\mathcal H'$. The graph $G'$ does not contain an independent set of size $h+1$, or equivalently, the complement of $G'$ is $K_{h+1}$-free. By Tur\'an's theorem \cite{Turan}, the number of edges of the complement of $G'$ is at most $(1-1/h)((h+1)s)^2/2$. Therefore, the number of edges of $G'$ is at least $\frac{(h+1)s((h+1)s-1)}{2}-\left(1-\frac{1}{h}\right)\frac{(h+1)^2s^2}{2}=\frac{(h+1)^2}{2h}\cdot\frac{n^{2/(h+1)}}{(\log n)^2}-\frac{(h+1)s}{2}$. Hence, $|B|\geq \pi(n)+\frac{h}{2}\cdot\frac{n^{2/(h+1)}}{(\log n)^2}$.

For proving the second statement our aim is to define a multiplicative basis of order $h$ for $[n]$ of the claimed size.
We are going to look for this basis in the form $B=P\cup X\cup Q$ where $P$ consists of the primes up to $n$, $X$ contains the integers up to $s^2$ and $Q$ contains certain 2-factor products of primes:
$$P=\{p\ |\ p\leq n\text{ and }p\text{ is a prime}\}, X=\{x\ |\ x\leq s^2\}, Q=\bigcup\limits_{-4\leq i\leq v}Q_i,$$
where the $Q_i$ sets (and $v$) are defined as follows. At first we are going to define $Q$ in the case $h\geq 14$.
Let $Q_{-1}=\{q_1q_2\ |\ q_1,q_2\in P,\ q_1\leq (h+1)^{-2}n^{1/(h+1)}, q_2\leq 2n^{1/(h+1)}   \}$ and $Q_{-2}=\{pq\ |\ p,q\in P, p\leq n/q^h, q\geq 2n^{1/(h+1)}       \} $. For defining $Q_{-3}$, let us divide the set $S$ of primes not greater than $2^{1.8}n^{1/(h+1)}$ into $r=\lfloor 0.61(h+1) \rfloor$ almost equal parts:  $S_{1},\dots,S_{r}.$ That is, for every $1\leq l\leq r$ we have $|S_l|=\left\lfloor \frac{\pi(2^{1.8}n^{1/(h+1)})}{r}      \right\rfloor$ or $\left \lceil \frac{\pi(2^{1.8}n^{1/(h+1)})}{r}      \right\rceil$, and $S$ is the disjoint union of the sets $S_1,\dots, S_r$. Let $Q_{-3}=S_1^2\cup\dots \cup S_r^2$. (For $h\geq 14$ let $Q_{-4}=\emptyset$.)

Let $v=\lfloor \log_2 (h+1) +\log_2 0.07 \rfloor$. Now, if $0\leq i\leq v$ let us divide the set $R_i$ of primes not greater than $2^{-i}n^{1/(h+1)}$ into $r_i=\lfloor 0.07\cdot2^{-i}(h+1) \rfloor$ almost equal parts: $R_{i,1},\dots,R_{i,r_i}$. 
That is, for every $1\leq l\leq r_i$ we have $|R_{i,l}|=\left\lfloor    \frac{\pi(2^{-i}n^{1/(h+1)})}{r_i}         \right\rfloor$ or $|R_{i,l}|=\left\lceil    \frac{\pi(2^{-i}n^{1/(h+1)})}{r_i}         \right\rceil$ and $R_i$ is the disjoint union of the sets $R_{i,1},\dots,R_{i,r_i}$. Let $Q_i=R_{i,1}^2\cup R_{i,2}^2\cup\dots \cup R_{i,r_i}^2$.

If $2\leq h\leq 13$, then let $Q=Q_{-2}\cup Q_{-4}$, where $Q_{-4}=\{pq\ |\ p,q\in P, p\leq n^{1/(h+1)}, q\leq 2n^{1/(h+1)}       \}$.

Now, we prove that $B$ is a multiplicative basis of order $h$ for $[n]$.
Let $a\leq n$ be arbitrary. Let us write $a$ as $a=p_1p_2\dots p_t$, where $p_1\geq p_2\geq\dots \geq p_t$ are the prime factors in the canonical form of $a$. At first we show that  $a\in (P\cup X)^h$  unless $h<t$ and $p_hp_{h+1}> s^2$. If $t\leq h$, then $a\in P^t\subseteq (P\cup X)^h$ trivially holds, so assume that $h<t$ and $p_hp_{h+1}\leq s^2$. Our aim is to distribute the primes appearing in the canonical form of $a$ into $h$ groups in such a way that in each group containing at least two elements the product of the primes is at most $s^2$.
The primes are going to be distributed into $h$ sets with a greedy algorithm. Let the products in these $h$ sets be $A_1,A_2,\dots,A_h$. At the beginning $A_1^{(0)}=A_2^{(0)}=\dots=A_h^{(0)}=1$. Then we put $p_1$ in the first set: $A_1^{(1)}:=p_1$. If $p_1,p_2,\dots,p_{l-1}$ are already distributed, then we put $p_l$ into the $j$-th group, if $A_j^{(l-1)}=\min \left(A_1^{(l-1)},A_2^{(l-1)},\dots,A_h^{(l-1)}\right)$, that is, if $A_j$ is currently one of the smallest products. (If there are more than one such $j$-s, we choose one arbitrarily.)
So, after the first $h$ steps we have $h$ many 1-factor products: $A_1^{(h)}=p_1,A_2^{(h)}=p_2,\dots,A_h^{(h)}=p_h$, then $p_{h+1}$ goes to the $h$-th group: $A_h^{(h+1)}=p_hp_{h+1}\leq s^2$.
We claim that by following this process at the end all of the products $A_1^{(t)},A_2^{(t)},\dots,A_h^{(t)}$ lie in $P\cup X$. For the sake of contradiction assume that at least one of them is not in $P\cup X$. Let $p_l=q$ be the first prime which created a product (with at least two prime factors) larger than $s^2$. Let us assume that after distributing the primes $p_1,p_2,\dots,p_{l-1}$ the products are $A_h\leq A_{h-1}\leq \dots \leq A_1$. Note that according to the indirect assumption $p_hp_{h+1}\leq s^2$ the number $l$ has to be at least $h+2$.
As $A_hq>s^2$, we have $s^2/q<A_h\leq A_{h-1}\leq\dots \leq A_1$. Hence, $(s^2/q)^h q<A_1A_2\dots A_hq\leq n$, thus
$$q>\left(\frac{s^{2h}}{n}\right)^{1/(h-1)}=\frac{n^{1/(h+1)}}{(\log n)^{2h/(h-1)}}.$$
Since $l\geq h+2$, we have $q\leq n^{1/(h+2)}$ which implies that $n^{1/(h+1)(h+2)}<(\log n)^{2h/(h-1)}$, however this contradicts the assumption $h\leq \sqrt{\frac{\log n}{12\log\log n}}$.

It is obtained that if $a\notin (P\cup X)^h$, then $p_hp_{h+1}>s^2$. Therefore, $p_1p_2\dots p_{h-1}<n/s^2$, which implies that $p_h\leq p_{h-1}< (n/s^2)^{1/(h-1)}=n^{1/(h+1)}(\log n)^{2/(h-1)}$. Therefore, $$p_{h+1}>s^2/p_{h}\geq n^{1/(h+1)} (\log n)^{-2h/(h-1)}$$
 and
$$p_1\leq n/p_{h+1}^h\leq n^{1/(h+1)} (\log n)^{2h^2/(h-1)}.$$
Summarizing these bounds we obtain that
\begin{multline}\label{eqprim}
n^{1/(h+1)} (\log n)^{2h^2/(h-1)}  \geq p_1\geq p_2\geq \dots \geq p_{h+1}\geq n^{1/(h+1)}(\log n)^{-2h/(h-1)}.
\end{multline}
Furthermore,
\begin{equation}\label{a'small}
a'\leq n/(p_1p_2\dots p_{h+1}) \leq n/(p_{h}p_{h+1})^{(h+1)/2}\leq n/s^{h+1}=(\log n)^{h+1},
\end{equation}
since the geometric mean of the numbers $p_1,\dots,p_{h+1}$ is bounded from below by the geometric mean of the two smallest elements: $p_h$ and $p_{h+1}$.

Hence, if $a\in [n]$, but $a\notin (P\cup X)^h$, then $a=p_1\dots p_{h+1}a'$, where the primes $p_1,\dots,p_{h+1}$ satisfy \eqref{eqprim} and $a'$ satisfies \eqref{a'small}.
We claim that if for all primes $p_1,\dots,p_{h+1}$ satisfying \eqref{eqprim} and $p_1p_2\dots p_{h+1}\leq n$ there exist some indices $1\leq i< j\leq h+1$ such that $p_ip_j\in Q$, then $B=P\cup X\cup Q$ is a multiplicative basis of order $h$ for $[n]$. To prove this, let us assume that $a=p_1\dots p_{h+1}a'$ satisfies these conditions and $p_ip_j\in Q$ for some $1\leq i< j\leq h+1$. Let $l\leq h+1$ be maximal such that $l\notin \{i,j\}$. Then $l\in \{h-1,h,h+1\}$, hence, $p_l\leq p_{h-1}\leq n^{1/(h+1)}(\log n)^{2/(h-1)}$. As $h\leq \sqrt{\frac{\log n}{12\log\log n}}$, $p_la'\leq n^{1/(h+1)}(\log n)^{(h^2+1)/(h-1)}<s^2$.
Let $q_1,\dots,q_{h-2}$ be the list of primes from $p_1,\dots,p_{h+1}$ excluding $p_i,p_j,p_l$ (only one appearance of each of them is excluded). Then $q_1,\dots,q_{h-2}\in P, p_ip_j\in Q, p_la'\in X$, so $a=q_1\dots q_{h-2}(p_ip_j)(p_la')\in B^h$.

It only  remains to show that for every primes $p_1,\dots,p_{h+1}$ satisfying \eqref{eqprim} and $p_1\dots p_{h+1}\leq n$ there exist some indices $1\leq i< j \leq h+1$ such that $p_ip_j\in Q$.

We start with the case $14\leq h$.
At first let us assume that $p_{h+1}\leq (h+1)^{-2}n^{1/(h+1)}$. If $p_h\leq 2n^{1/(h+1)}$, then $p_{h+1}p_h\in Q_{-1}$, and we are done. Otherwise, $p_h> 2n^{1/(h+1)}$ and $p_{h+1}\leq n/p_h^h$, hence, $p_hp_{h+1}\in Q_{-2}$. Thus it can be assumed that $p_{h+1}>(h+1)^{-2}n^{1/(h+1)}$.

Let us denote the multiset of $p_1,\dots,p_{h+1}$ by $T$. For $i\geq 0$ let $N_i$ denote the number of such elements of $T$ that are at most $2^{-i}n^{1/(h+1)}$. At first let us assume that there exists some $0\leq i\leq v$ such that $N_i>0.07\cdot 2^{-i} (h+1)\geq r_i$. Since $T$ contains more than $r_i$ elements of the set $R_i$, by the pigeonhole principle there exist some indices $l_1$ and $l_2$ such that $p_{l_1},p_{l_2}\in R_{i,j}$ for some $j$. Then $p_{l_1}p_{l_2}\in Q_i$, and we are done.

Now let us assume that for every $0\leq i\leq v$ we have $N_i\leq  0.07\cdot 2^{-i} (h+1)$. Specially, $N_v\leq 1$, that is, $T$ contains at most one element (namely, $p_{h+1}$) less than $2^{-v}n^{1/(h+1)}$, however, this element is at least $(h+1)^{-2}n^{1/(h+1)}$.
Let the multiset $T_1$ contain those elements of $T$ that are at most $n^{1/(h+1)}$, the remaining elements of $T$ are in $T_2$. Note that $h+1=|T|=|T_1|+|T_2|$.

 Now, a lower bound is going to be given for $\prod\limits_{p_i\in T_1} p_i$. Since all the elements of $T_1$ except $p_{h+1}$ are in the interval $(2^{-v}n^{1/(h+1)},n^{1/(h+1)}]$, the double-counting of the size of the set
$$\{(i,j)\ |\ p_i\leq 2^{-j}n^{1/(h+1)},p_i\in T_1\setminus\{p_{h+1}\},0\leq j\text{ is an integer}\}$$
yields the estimate
$$\prod\limits_{p_i\in T_1\setminus \{p_{h+1}\}} p_i\geq n^{(|T_1|-1)/(h+1)}2^{-\sum\limits_{i=0}^v N_i}.$$
Therefore,
\begin{multline*}
\prod\limits_{p_i\in T_1} p_i\geq n^{(|T_1|-1)/(h+1)}2^{-\sum\limits_{i=0}^v N_i}p_{h+1}\geq n^{(|T_1|-1)/(h+1)}2^{-\sum\limits_{i=0}^v 0.07\cdot 2^{-i} (h+1)}p_{h+1} \geq \\
\geq n^{(|T_1|-1)/(h+1)}2^{-0.14(h+1)}p_{h+1}\geq n^{|T_1|/(h+1)}2^{-0.69(h+1)}, 
\end{multline*}
where we used that $(h+1)^2\leq 2^{0.55(h+1)}$ for every $h\geq 14$. Note that   $|T_1|=N_0\leq 0.07(h+1)$. As $p_1\dots p_{h+1}\leq n$, the following upper bound is obtained for the product of the elements of $T_2$:
$$\prod\limits_{p_i\in T_2}p_i\leq n^{|T_2|/(h+1)}2^{0.69(h+1)}.$$
Therefore, $T_2$ contains less than $0.39(h+1)$ elements larger than $2^{1.8}n^{1/(h+1)}$. Hence, more than $0.61(h+1)$ elements of $T_2$ are at most $2^{1.8}n^{1/(h+1)}$. Then, by the pigeonhole principle two elements of $T_2$ lie in the same set $S_j$, therefore their product is in $Q_{-3}$ and we are done.

Finally, if $2\leq h\leq 13$, then $p_h\geq 2n^{1/(h+1)}$ implies $p_hp_{h+1}\in Q_{-2}$ and $p_h\leq 2n^{1/(h+1)}$ implies $p_hp_{h+1}\in Q_{-4}$.

Hence, it is shown that $B$ is a multiplicative basis of order $h$ for $[n]$.

Finally, an upper bound will be given for the size of $B$. Clearly, $|P|=\pi(n)$, $|X|\leq s^2$.

For the size of $Q_{-1}$ we have that $|Q_{-1}|\leq 1.26^2\cdot 2s^2\leq 0.3 hs^2$ for every $h\geq 14$.


As $Q_{-2}=\{pq\ |\ p,q\in P, p\leq n/q^h, q\geq 2n^{1/(h+1)}       \}=\bigcup\limits_{1\leq j} \{pq\ |\ p,q\in P, p\leq n/q^h, 2^jn^{1/(h+1)}\leq q<2^{j+1}n^{1/(h+1)}       \}\subseteq \bigcup\limits_{1\leq j} \{pq\ |\ p,q\in P, p\leq 2^{-jh}n^{1/(h+1)},q\leq 2^{j+1}n^{1/h+1} \}$, we have that $|Q_{-2}|\leq 1.26^2\sum\limits_{1\leq j} 2^{-jh+j+1}(h+1)^2s^2=1.26^2\frac{(h+1)^22^{2-h}}{1-2^{1-h}}s^2\leq 14.3hs^2$ for every $h\geq 2$.

If $h\geq 14$, then $9\leq r$, so $\lfloor 0.61 (h+1) \rfloor \geq 0.61(h+1)(9/10)$. Hence, $|Q_{-3}|\leq 1.26^2 (5/4)\frac{(h+1)^22^{2\cdot 1.8}}{0.61(h+1)} s^2$, that is, we have $|Q_{-3}|=\frac{{1.26^2}(10/9){2^{3.6}}}{0.61}\cdot\frac{15}{14}hs^2\leq 37.6 hs^2$.


If $0\leq i\leq v$, then $r_i\geq 1$, so $\lfloor 0.07\cdot 2^{-i}(h+1)\rfloor \geq 0.07\cdot 2^{-i}(h+1)/2$. Therefore, $|Q_i|\leq 2\frac{1.26^2}{0.07}2^{-i}(h+1)s^2$, so for the size of the union of the sets $Q_0,\dots,Q_v$ we obtain that:
$\sum\limits_{0\leq i\leq v} |Q_i|\leq 97.2hs^2$ for every $h\geq 14$.

If $2 \leq h \leq 13$, then $|Q_{-4}|\leq 1.26^2 2(h+1)^2s\leq 47.9 hs^2$.

Hence, $|B|\leq |P|+|X|+|Q|=\pi(n)+150.4hs^2$, if $14\geq h$ and $|B|\leq |P|+|X|+|Q|\leq \pi(n)+63.2 hs^2$, if $2\leq h\leq 13$. $\blacksquare $


Now we continue with the problem of Erd\H{o}s, estimating $F_h(n)$. We start with proving two lemmas.

\begin{lemma}\label{tinter}
Let $k$ be a fixed positive integer. Let $S$ be a set of size $n\ge 2k^2$. Then for every $1\le t<k$ one can choose $l$ many $k$-element subsets $S_1,S_2,\dots ,S_l\subset S$ such that for every $i\not =j$ we have $|S_i\cap S_j|<t$ and $l\ge \left( \frac{n}{2k}\right) ^t$.
\end{lemma}

\noindent
{\bf Proof of Lemma \ref{tinter}}. By Bertrand's postulate there exists a prime number $q$ between $\frac{n}{2k}$ and $\frac{n}{k}$. It can be supposed that $S\supseteq\mathbb{F}_q\times [k]$. That is, it can be assumed that $S$ contains $k$ disjoint copies of $\mathbb{F}_q$. All the $k$-element sets are going to contain one element from each copy of $\mathbb{F}_q$ in such a way that the intersection of any two of them has size smaller that $t$.
 These $q^t\geq \left( \frac{n}{2k}\right) ^t $ suitable sets $S_i$ are defined in the following way: Let $p(x)=a_0+a_1x+\dots +a_{t-1}x^{t-1}$, where $a_0,a_1,\dots,a_{t-1}\in [0,q-1)$.
$$S_{p(x)}:=\bigcup\limits_{1\leq i\leq k}(p(i),i).$$
It remains to prove that for different polynomials $p_1(x)$ and $p_2(x)$ we have $|S_{p_1(x)}\cap S_{p_2(x)}|<t $. For the sake of contradiction, let us assume that $|S_{p_1(x)}\cap S_{p_2(x)}|\geq t $. Then there exist
 $1\le x_1<x_2<\dots <x_t\le k$ such that $p_1(x_i)=p_2(x_i)$ for every $1\le i\le t$, which contradicts that the degree of $p_1-p_2$ is at most $t-1$. $\blacksquare $

\begin{lemma}\label{inj}
Let $A\subseteq [n]$ possessing property $\mathcal{P}_h$ and  $B\subset [n]$.
Then there exists a one-to-one mapping $A\cap B^h\to B$ such that for $a\to b$ there exist integers $b_2,\dots ,b_h\in B$ such that $a=bb_2\dots b_h$. As a special case, if $B$ is a multiplicative basis of order $h$ for $[n]$, then there is a one-to-one mapping $A\to B$ such that for $a\to b$ there exist integers $b_2,\dots, b_h\in B$ such that $a=bb_2\dots b_h$.
\end{lemma}
\noindent


\noindent
{\bf Proof of Lemma \ref{inj}}. Let us write each element in $A_0=A\cap B^h$ as a product of $h$ (not necessarily distinct) elements of $B$. (If there are more than one possibilities, let us choose one arbitrarily.) Let $a\in A$, and the representation of $a$ be $a=b_1^{\lambda_1}\dots b_k^{\lambda_k}$ where $\lambda_1+\dots+\lambda_k=h$. We claim that for some $1\leq i\leq k$ the factor $b_i$ appears in the representation of any element of $A_0\setminus \{a\}$ at most $\lambda_i-1$ times. For the sake of contradiction assume that for every $1\leq i\leq k$ there is an $a_i\in A_0\setminus \{a\}$ such that $b_i$ appears in the representation of $a_i$ at least $\lambda_i$ times. Let $a_1',\dots,a_l'$ be the distinct elements of the multiset $\{a_1,\dots,a_k\}$. (That is, the elements are listed without repetition, $l\leq k$.)
Then $a|a_1'\dots a_l'$, which contradicts that $A$ possesses property $\mathcal{P}_h$, since $l\leq k\leq h$. Therefore, there is an $i$ for which the multiplicity of $b_i$ in the representation of $a$ is maximal. Let us assign such a $b_i$ to $a$. Clearly, this is a one-to-one mapping.

In the special case when $B$ is a multiplicative basis of order $h$ for $[n]$, we have $A_0=A\cap B^h=A$. $\blacksquare $


\noindent
Now, we are ready to prove Theorem~\ref{fhveges}.
\smallskip

\noindent
{\bf Proof of Theorem \ref{fhveges}}. Let $n^{\frac{1}{h+1}}(\log n)^{-1}=s$.
At first we prove the lower bound. Let $S$ be the set of primes not greater than  $n^{1/(h+1)}$. Since $|A|\geq 2h^2$, Lemma~\ref{tinter} implies that we can choose $\left(\frac{|S|}{2(h+1)}\right)^2$ many subsets of $S$ of size $h+1$ in such a way that the intersection of any two of them contains at most one element. Let these subsets be $S_1,\dots,S_m$, where $m\geq \left(\frac{|S|}{2(h+1)}\right)^2$. Now, let $s_i=\prod\limits_{s\in S_i}s$ for every $1\leq i\leq m$ and $A=\{s_i:\ 1\leq i\leq m\}\cup \{q\ |\ n^{1/(h+1)}<q\leq n,q\text { is a prime}\}$. We claim that $A$ possesses property $\mathcal{P}_h$. Since, if $a,a_1,\dots,a_h$ are distinct elements of $A$, and $a$ is a product of $h+1$ primes, then every $a_i$ is divisible by at most one of these prime factors implying that $a$ can not divide $a_1a_2\dots a_h$. On the other hand, if $a\in A$ is a prime, then $a>n^{1/(h+1)}$ and there is no other element in $A$ which is divisible by $a$, hence $a\nmid a_1a_2\dots a_h$.
Furthermore $|A|\geq \pi(n)-\pi(n^{1/(h+1)})+\left(\frac{\pi(n^{1/(h+1)})}{2(h+1)}\right)^2>\pi(n)+0.2s^2$.


Now we continue with the upper estimate. Let $A\subseteq [n]$ be a set possessing property $\mathcal{P}_h$. Lemma~\ref{inj} implies that $|A|\leq G_h(n)$. In the proof of Theorem~\ref{MBveges} we showed that for $h\leq 6$ we have $G_h(n)\leq \pi(n)+63.2hs^2$, therefore, $F_h(n)\leq \pi(n)+379.2s^2$ also holds. From now on, we assume that $7\leq h$.

Let $P$ be the set of the primes up to $n$ and $X$ contain the integers up to $s^2$:
$$P=\{p\ |\ p\leq n\text{ and }p\text{ is a prime}\}, X=\{x\ |\ x\leq s^2\}.$$

Now a  mapping from a subset of $A$ to $P\cup X$ is going to be defined in 3 steps:
\begin{itemize}
\item[(i)] If $a\in A$ and there exists a prime $p\in P(s^2)$ and an exponent $\alpha$ such that $p^\alpha|a$, but $p^\alpha\nmid a'$ for every $a\ne a'\in A$, then let us assign such a $p$ to $a$.
\item[(ii)] Let us write each element of $A\cap (P\cup X)^h$ as a product of $h$ elements from $P\cup X$. If $a\in A$ does not have an image yet, moreover, there exists a $y\in P\cup X$ and an $\alpha\in\mathbb{Z}^+$ such that $y$ occurs $\alpha$ times in the representation of $a$, but it occurs at most $\alpha-1$ times in the representation of any other $a'\in A\cap (P\cup X)^h$, then let us assign such a $y$ to $a$.
\item[(iii)] Finally, if an element $a\in A$ does not have an image yet, but there exists an $x\in X$ such that $x|a$, but $x\nmid a'$ for every $a\ne a'\in A$, then let us assign such an $x$ to $a$.
\end{itemize}
Let $A_1\subseteq A$ contain those elements of $a$ that has an image and $A_2:=A\setminus A_1$.
If an element of $P\cup X$ is assigned to more than one element of $A_1$, then it has to be a prime which is at most $s^2$, and it is assigned to exactly two elements: one according to rule (i) and one according to rule (ii). Therefore, $|A_1|\leq |P|+2|X|\leq \pi(n)+2s^2$.
According to Lemma~\ref{inj} we have $A\cap (P\cup X)^h\subseteq A_1$.

Finally, our aim to show that $|A_2|\leq357.2s^2$. Let $a\in A_2$. As we have seen in the proof of Theorem~\ref{MBveges}, since $a\notin (P\cup X)^h$, the number $a$ can be written as $a=p_1p_2\dots p_{h+1}a'$, where the primes $p_1,p_2,\dots,p_{h+1}$ satisfy the condition
 \begin{multline}
n^{1/(h+1)} (\log n)^{2h^2/(h-1)}  \geq p_1\geq p_2\geq \dots \geq p_{h+1}\geq n^{1/(h+1)}(\log n)^{-2h/(h-1)},
\end{multline}
moreover $a'\leq (\log n)^{h+1}$ and $p_1p_2\dots p_{h+1}\leq n$. Let us denote the multiset $\{p_1,\dots,p_{h+1}\}$ by $T=T_a$. Note  that all elements of $T_a$ are less than $s^2$. We claim that for every $a\in A_2$ the multiset $T_a$ contains $h+1$ distinct primes. For the sake of contradiction assume that the multiset $T_a$ contains $\lambda_1$ many $q_1$'s, $\lambda_2$ many $q_2$'s, and so on, $\lambda_t$ many $q_t$'s, where $q_1,q_2,\dots,q_t$ are distinct primes and $t\leq h$. That is, $a=q_1^{\lambda_1}\dots q_t^{\lambda_t}a'$, where $\lambda_1+\dots+\lambda_t=h+1$. As $a\notin A_1$, there exist $b_1,\dots,b_t,b_{t+1}\in A\setminus \{a\}$ such that $q_1^{\lambda_1}|b_1,\dots,q_t^{\lambda_t}|b_t,a'|b_{t+1}$. Let the multiset $\{b_1,\dots,b_{t+1}\}$ contain the pairwise different elements $c_1,\dots,c_u$, where $u\leq t+1$. Then $a|c_1\dots c_u$, since $q_1^{\lambda_1},\dots,q_t^{\lambda_t}$ and $a'$ are pairwise coprimes. If $t+1\leq h$, then this contradicts the assumption that $A$ possesses property $\mathcal{P}_h$. Therefore, it can be assumed that $t=h$. Then $a=q_1^2q_2\dots q_h a'$, and without the loss of generality, it can be assumed that $q_2\geq q_3\geq \dots\geq q_h$. Since $q_h\in\{p_{h-1},p_h,p_{h+1}\}$, we have $q_h\leq n^{1/(h+1)}(\log n)^{2/(h-1)}$, hence, $q_ha'\leq s^2$, that is, $q_ha'\in X$. As $a\notin A_1$, there exist $b_1,\dots,b_h\in A\setminus \{a\}$ such that $q_1^2|b_1,q_2|b_2,\dots,q_{h-1}|b_{h-1},q_ha'|b_h$.  Let the multiset $\{b_1,\dots,b_{h}\}$ contain the pairwise different elements $c_1,\dots,c_u$, where $u\leq h$. Then $a|c_1\dots c_u$, since $q_1^2, q_2,\dots,q_{h-1},q_ha'$ are pairwise coprimes, however this contradicts the assumption that $A$ possesses property $\mathcal{P}_h$. Therefore, for every $a\in A_2$ the multiset $T_a$ contains $h+1$ distinct primes.

Now we claim that for any two different elements  $a,b\in A_2$ the  intersection of $T_a=\{p_1,p_2,\dots, p_{h+1}\}$ and $T_b$ contains at most one prime, that is, $|T_a\cap T_b|\leq 1$. For the sake of contradiction assume that for some $a,b\in A_2$ we have $|T_a\cap T_b|\geq 2$. Namely, let $1\leq i<j\leq h+1$ be the two indices for which $p_i,p_j\in T_b$. Let $l$ be maximal such that $l\notin \{i,j\}$. Then $l\in \{h-1,h,h+1\}$, thus $p_la'\in X$. Let $\{q_1,q_2,\dots,q_{h-2}\}=T_a\setminus \{p_i,p_j,p_l\}$. As $a\in A_2$, for every $1\leq m\leq h-2$ there exists $b_m\in A$ such that $q_m|b_m$ and there exists $b_{h-1}\in A$ such that $p_la'|b_{h-1}$.
Let $c_1,\dots,c_u$ be the distinct elements of the multiset $\{b,b_1,\dots,b_{h-1}\}$, so $u\leq h$.
Then $a|c_1\dots c_u$, since $q_1,\dots,q_{h-2}, p_ip_j,p_la'$ are pairwise coprimes. This contradicts the assumption that $A$ possesses property $\mathcal{P}_h$.

Therefore,  each $T_a$ (where $a\in A_2$) contains $h+1$ distinct primes, moreover the intersection of $T_a$ and $T_b$ contains at most one element (if $a,b\in A_2$ and $a\ne b$).

Let $C$ contain those elements $a$ of $A_2$ for which $\min \{T_a\}<(h+1)^{-2}n^{1/(h+1)}$.
Let $Q_{-1}=\{q_1q_2\ |\ q_1,q_2\in P,\ q_1\leq (h+1)^{-2}n^{1/(h+1)}, q_2\leq 2n^{1/(h+1)}   \}$ and $Q_{-2}=\{pq\ |\ p,q\in P, p\leq n/q^h, q\geq 2n^{1/(h+1)}       \} $.
Let $a=p_1p_2\dots p_{h+1}a'\in C$.
If $p_h\leq 2n^{1/(h+1)}$, then $p_{h+1}p_h\in Q_{-1}$. Otherwise, $p_h> 2n^{1/(h+1)}$ and $p_{h+1}\leq n/p_h^h$, hence, $p_hp_{h+1}\in Q_{-2}$.
Therefore $|C|\leq |Q_{-1}|+|Q_{-2}|\leq31.8s^2$, where the upper bounds for the sizes of $Q_{-1}$ and $Q_{-2}$ can be obtained similarly as in the proof of Theorem~\ref{MBveges}.

From now on, it is assumed that $a\in A_2\setminus C$.
For $i\geq 0$ let us denote by $P_i$ the set of primes not greater than $2^{-i}n^{1/(h+1)}$.
Moreover, let $N_i=N_i(a)$ denote the size of $T_a\cap P_i$. 
Let $v=\left\lceil \frac{\log_2 (h+1)+\log_2 0.17}{1-\log_2 1.2}\right\rceil-1\geq 0$.
Let $A'_{i}$ be the set of those elements $a\in A_2$ for which $N_i(a)\geq r_i=0.17\cdot2^{-i}1.2^i(h+1)$.  
If $i\leq v$, then $r_i> 1$ and $\binom{N_i(a)}{2}\geq \binom{r_i}{2}>0$. Each 2-element subset of $P_i$ is contained in at most one $T_a$. However each $T_a$ contains at least $\binom{r_i}{2}$ many 2-element subsets of $P_i$, therefore
$$|A'_i|\leq \frac{\binom{|P_i|}{2}}{\binom{r_i}{2}} =\frac{\binom{\pi(2^{-i}n^{1/(h+1)})}{2}}{\binom{r_i}{2}}\leq 1.26^2\cdot 0.17^{-2} 1.2^{-2i}s^2.$$
Furthermore, $\sum |A_i'|\leq \frac{1.26^2\cdot 0.17^{-2}}{1-1.2^{-2}}s^2\leq 179.8s^2$.

Now, if $a\in A^*=A_2\setminus \bigcup\limits_{0\leq i\leq v} A'_i$, then in $T_a$ the number of elements smaller than $2^{-i}n^{1/(h+1)}$ is $N_i\leq r_i=0.17\cdot2^{-i}1.2^i(h+1)$ for every $0\leq i\leq v$.
Specially, $N_v\leq 1$, that is, $p_1,p_2,\dots,p_h$ are all at least $2^{-v}n^{1/(h+1)}$.
Let $T_{a}^{(1)}$ contain those elements of $T_a$ that are at most $n^{1/(h+1)}$ and let $T_{a}^{(2)}=T_a\setminus T_{a}^{(1)}$.


Now, a lower bound is going to be given for $\prod\limits_{p_i\in T_{a}^{(1)}} p_i$. Since all elements of $T_{a}^{(1)}$ (possibly) except $p_{h+1}$ are in the interval $(2^{-v}n^{1/(h+1)},n^{1/(h+1)}]$, the double-counting of the size of the set
$$\{(i,j)\ |\ p_i\leq 2^{-j}n^{1/(h+1)},p_i\in T_{a}^{(1)}\setminus\{p_{h+1}\},0\leq j\text{ is an integer}\}$$
yields the estimate
$$\prod\limits_{p_i\in T_{a}^{(1)}\setminus \{p_{h+1}\}} p_i\geq n^{(|T_1|-1)/(h+1)}2^{-\sum\limits_{i=0}^v N_i}.$$
Therefore,
\begin{multline*}
\prod\limits_{p_i\in T_{a}^{(1)}} p_i\geq n^{(|T_{a}^{(1)}|-1)/(h+1)}2^{-\sum\limits_{i=0}^v N_i}p_{h+1}\geq n^{(|T_{a}^{(1)}|-1)/(h+1)}2^{-\sum\limits_{i=0}^v 0.17\cdot 2^{-i}1.2^i (h+1)}p_{h+1} \geq \\
\geq n^{(|T_{a}^{(1)}|-1)/(h+1)}2^{-c(h+1)}p_{h+1}\geq n^{|T_a^{(1)}|/(h+1)}2^{-(c+0.75)(h+1)}, 
\end{multline*}
where $c=0.17\cdot (1-1.2/2)^{-1}=0.425$ and we used that $(h+1)^2\leq 2^{0.75(h+1)}$ for every $h\geq 7$. Note that   $|T_a^{(1)}|=N_0\leq 0.17(h+1)$. As $p_1\dots p_{h+1}\leq n$, the following upper bound is obtained for the product of the elements of $T_a^{(2)}$:
$$\prod\limits_{p_i\in T_a^{(2)}}p_i\leq n^{|T_a^{(2)}|/(h+1)}2^{1.175(h+1)}.$$
Therefore, $T_2$ contains at most $0.51(h+1)$ elements larger than $2^{2.3}n^{1/(h+1)}$. 
So at least $0.49(h+1)$ elements of $T_2$ are at most $2^{2.3}n^{1/(h+1)}$.
Hence, $|A^*|\leq \frac{    \binom{\pi(2^{2.3}n^{1/(h+1)})}{2}            }{     \binom{0.49(h+1)}{2}             } \leq\frac{1.26^2\cdot 2^{2\cdot 2.23}}{0.49^2}s^2\leq 145.6 s^2$.

Therefore, it is obtained that $|A|\leq \pi(n)+359.2s^2$, if $h\geq 7$. $\blacksquare $

\noindent
We note that in the proofs of Theorem~\ref{MBveges} and Theorem~\ref{fhveges} with a more careful and lengthier calculation better constants can be obtained, especially, if $h$ is large enough.


\section{Infinite case}

In this section the following lemma of Erd\H os is going to be used (\cite{ep}).

\begin{lemma}\label{erdos}
The set $B=\left\{k: k\leq n^{2/3}\right\} \cup \{ p: p\leq n\text{ and $p$ is a prime}\} $ forms a multiplicative basis of order 2 for $[n]$.
\end{lemma}
\noindent
Our first lemma generalizes Erd\H os' previously mentioned lemma.
\begin{lemma}\label{basis}
Let $h\geq 2$. The set $B^{(h)}=\{k: k\leq n^{\frac{2}{h+1}}\} \cup \{ p: p\leq n\text{ and $p$ is a prime}\} $ forms a multiplicative basis of order $h$ for $[n]$.
\end{lemma}

\noindent
{\bf Proof of Lemma \ref{basis}}.
We prove the statement by induction on $h$. The base case $h=2$ was shown by Erd\H os.

Now let us suppose that for every $N$ the set $\{k: k\leq N^{2/h}\} \cup \{ p: p\leq N\text{ and $p$ is a prime}\} $ forms a multiplicative basis of order $h-1$ for $[N]$. We show that $B^{(h)}=\{k: k\leq n^{\frac{2}{h+1}}\} \cup \{ p: p\leq n\text{ and $p$ is a prime}\} $ forms a multiplicative basis of order $h$ for $[n]$. Let $m\leq n$. If there exists a prime divisor $p$ of $m$ such that $p>n^{\frac{1}{h+1}}$, then $m=p\cdot\frac{m}{p}$, where $\frac{m}{p}\leq n^{\frac{h}{h+1}}$. Therefore, using the induction step for $N=n^{\frac{h}{h+1}}$ we get that $\frac{m}{p}=b_2\dots b_h$ such that either $b_i\leq N^{\frac{2}{h}}=n^{\frac{2}{h+1}}$ or $b_i$ is a prime, so $m=b_1b_2\dots b_h$ for some $b_i\in B^{(h)}$.

If every prime divisor of  $m$ is at most $n^{\frac{1}{h+1}}$, then let $m=p_1p_2\dots p_s$ such that $p_1\geq p_2\geq \dots \geq p_s$. We show that the multiset $\{ p_1,p_2,\dots ,p_s\} $ can be split into $h$ parts, $A_1\cup A_2 \cup \dots \cup A_h$, such that every number of the form $b_i=\prod\limits_{p\in A_i}p$ is at most $n^{\frac{2}{h+1}}$. 
Let $A_1^{(h)}=\{p_1\},\dots ,A_h^{(h)}=\{p_h\}$. Now assume that for some $h\leq i<s$ we have already defined the multisets $A_1^{(i)},\dots ,A_h^{(i)}$. Let $b_g^{(i)}=\prod\limits_{p\in A_g^{(i)}}p$ and $j$ is chosen in such a way that $\min \{ b_1^{(i)},b_2^{(i)},\dots ,b_h^{(i)}\}=b_j^{(i)}$. Then let $A_g^{(i+1)}=A_g^{(i)}$ for every $g\not =j$ and  $A_j^{(i+1)}=A_j^{(i)}\cup \{ p_{i+1}\} $. We claim that $b_{j+1}^{(i+1)}\leq n^{\frac{2}{h+1}}$. For the sake of contradiction let us  assume that $b_{j+1}^{(i+1)}> n^{\frac{2}{h+1}}$. Then $n\geq m\geq b_1^{(i+1)}\dots b_{h}^{(i+1)}> b_1^{(i+1)}\dots b_{j-1}^{(i+1)}b_{j+1}^{(i+1)}\dots b_h^{(i+1)}n^{\frac{2}{h+1}}$, therefore $n^{\frac{h-1}{h+1}}>b_1^{(i)}\dots b_{j-1}^{(i)}b_{j+1}^{(i)}\dots b_h^{(i)}$. Thus $\min \{ b_1^{(i)},b_2^{(i)},\dots ,b_h^{(i)}\}<n^{\frac{1}{h+1}}$. Hence $b_j^{(i+1)}=b_j^{(i)}p_{i+1}<n^{\frac{2}{h+1}}$ is a contradiction. Thus always adding the following prime to the set in which the product is currently the smallest gives us an appropriate representation. $\blacksquare $

\smallskip

\noindent
Now, we prove Theorem~\ref{MBliminf}.

\noindent
{\bf Proof of Theorem \ref{MBliminf}}.
We start with proving the first statement by induction on $h$ for every $h\geq 1$. First of all, note that the unique multiplicative basis of order 1 for $\mathbb{Z}^+$ is $B=\mathbb{Z}^+$, hence, for $h=1$ the statement is trivially true. Now assume that $h\geq 2$ and for $h-1$ the statement holds.
Let  $B\subseteq \mathbb{Z}^+$ be a multiplicative basis of order $h$ for $\mathbb{Z}^+$. Without the loss of  generality it can be assumed that $B$ is not a multiplicative basis of order $h-1$, otherwise the statement follows from the induction hypothesis. So, it can be supposed that there exists some $m\in \mathbb{Z}^+\setminus B^{h-1}$.
Clearly, all the primes (and 1) have to belong to $B$.
Now let $n$ be an arbitrary integer large enough. Let $m<p\leq n/m$ be an arbitrary prime. Since $pm\in B^h$, the number $pm$ can be written as $pm=b_1b_2\dots b_h$ in such a way that $b_1,b_2,\dots,b_h\in B$. As $p$ is a prime, it divides some $b_i$, so let us assume that $p|b_1$. Then $b_1>p$, since $b_1=p$ would imply that $m=b_2b_3\dots b_h\in B^{h-1}$. Therefore, $b_1\in B$ is a multiple of $p$, moreover, $b_1/p\leq m$.
Hence,  $b_1$ is a composite number and has a unique prime factor larger than $m$. For each prime from the interval $(m,n/m)$ we get such an element of $B$ and these elements are distinct, thus
$B(n)\geq \pi(n) +\pi(n/m)-\pi(m)$. Hence, $\liminf \frac{|B(n)|}{{\frac{n}{\log n}}}\geq 1+1/m$.

To prove the second statement, it is enough to do so in the special case $h=2$, since a multiplicative basis of order 2 is a multiplicative basis of order $h$ for every $h\geq 2$. Let $\varepsilon>0$ be arbitrary.
We are going to find an increasing sequence $(n_i)_{i=1}^\infty$ of positive integers and sets $B_i\subseteq [n_i]$ in such a way that the following conditions hold for every $i\geq 1$:
\begin{itemize}
\item[(i)] $B_i$ is a multiplicative basis for $[n_i]$,
\item[(ii)] $|B_i|<(1+\varepsilon)\frac{n_i}{\log n_i}$,
\item[(iii)] $B_{i}\cap[n_{i-1}]=B_{i-1}$.
\end{itemize}
If such numbers and sets are found, then let us define  a sequence of positive integers by $ B:=\displaystyle\bigcup\limits_{i=1}^\infty B_i$.
We claim that $B$ is a multiplicative basis for $\mathbb{Z}^+$ satisfying that $\liminf \frac{|B(n)|}{{\frac{n}{\log n}}}<1+\varepsilon$.
At first we show that $B$ is a multiplicative basis. Let $a\in\mathbb{Z}^+$ be arbitrary. If $i$ is large enough, then $a\in [n_i]$. Since $B_i$ is a multiplicative basis for $[n_i]$, there exist $b,c\in B_i$ such that $a=bc$. As $B_i\subseteq B$, the number $a$ is a product of two elements of $B$. This is true for every $a\in\mathbb{Z}^+$, so $B$ is a multiplicative basis.

Condition (iii) implies that $B(n_i)=B_i$, hence, by condition (ii) it follows that for every $i\geq 1$ we have
$$\frac{|B(n_i)|}{\frac{n_i}{\log n_i}}<1+\varepsilon.$$
From this the desired statement follows.

Now it remains to find appropriate $n_i$ numbers and $B_i$ sets.
Let $n_0=N=\lceil \max((32/(3\varepsilon))^2, (16/\varepsilon)^2)\rceil$ and $B_0=[n_0]$. Now we define the numbers $n_i$ and the sets $B_i$ (for $i\geq 1$) satisfying conditions (i), (ii), (iii) recursively. Let us assume that $n_i$ and $B_i$ are already defined in such a way that $B_i\subseteq [n_i]$ is a multiplicative basis for $[n_i]$. Our aim is to find $n_{i+1}>n_i$ and $B_{i+1}\subseteq [n_{i+1}]$ satisfying conditions (i), (ii), (iii). For simplicity let us introduce the notion $x:=n_i,y:=n_{i+1}$.
Let us define $B_{i+1}$ in the following way:
\begin{multline*}
B_{i+1}=B_i\cup\{i\ |\ x<i\leq y^{2/3}x\}\cup\{pv\ |\ y^{2/3}<p\leq y/x,p \text{ is a prime, }v\leq\sqrt{x}\}\cup \\
\cup \{pv\ |\ y/x<p\leq y/N,p\text{ is a prime}, v\leq\sqrt{ y/p}\}\cup \\
\cup \{p\ |\ y/N<p\leq y,p\text{ is a prime}\}.
\end{multline*}
If $y>x^2$, we have $\min(y^{2/3},y/x)>x$, so every element of $B_{i+1}\setminus B_i$ is larger than $x$, therefore condition (iii) holds.

Now we show that $B_{i+1}$ is a multiplicative basis for $[y]$. Let $a\leq y$ be arbitrary. According to Lemma~\ref{erdos} the number $a$ can be represented in the form $a=uv$, where $v\leq u$ and either $u\leq y^{2/3}$, or $u>y^{2/3}$ is a prime. At first assume that $u\leq y^{2/3}$. If $x<v$, then both $u$ and $v$ lie in $(x,y^{2/3}x]$, so $u,v\in B_{i+1}$ and $a=uv\in B_{i+1}^2$. If $v\leq x$, then we distinguish two cases.
\begin{enumerate}
\item If $x<uv$, then $a=1\cdot (uv)$ is a good representation, since $uv$ lies in $(x,y^{2/3}x]$.
\item If $uv\leq x$, then $a=uv$ can be written as a product of two elements from the set $B_i\subseteq B_{i+1}$, since $B_i$ is a multiplicative basis for $[x]$ by the induction hypothesis.
\end{enumerate}

Secondly let us assume that $u>y^{2/3}$ is a prime, denote it by $p$. As the first case let $y^{2/3}<p\leq y/x $. Since $a\leq y$, we have that $v=a/p\leq y/p\leq y^{1/3}$. If $x<v$, then $v\in B_{i+1}$, so $a=pv\in B_{i+1}^2$. If $v\leq x$, then $v=v_1v_2$ for some $v_1,v_2\in B_i$, since $B_i$ is a multiplicative basis for $[x]$. Without the loss of generality it can be assumed that $v_1\leq v_2$. Then $v_1\leq \sqrt{v_1v_2}=\sqrt{v}\leq \sqrt{x}$, therefore both  $pv_1$ and $v_2$ lies in $B_{i+1}$, hence $a=(pv_1)\cdot v_2\in B_{i+1}^2$.

Now, as the second case let $y/x<p$. If $y/N \leq p$, then $v=a/p\leq y/p\leq N$, so $a=p\cdot v$ is a good representation, since $[N]\subseteq B_{i+1}$ and $p\in B_{i+1}$. Finally, if $y/x<p<y/N$, then 
$v=y/p< x$. Since $B_{i}$ is a multiplicative basis for $[x]$, there exist some $v_1,v_2\in B_i$ such that $v=v_1v_2$. It can be assumed that $v_1\leq v_2$ and in this case $v_1\leq \sqrt{v}\leq \sqrt{y/p}$. Therefore, $a=(pv_1)\cdot v_2\in B_{i+1}^2$. Thus we obtained that condition (i) holds.

Finally, it is going to be proved that $B_{i+1}$ and $n_{i+1}$ satisfies condition (ii), as well.
If $x^4<y$, then
$$|B_i\cup\{i\ |\ x<i\leq y^{2/3}x\}|\leq y^{2/3}x<y^{11/12}<\frac{\varepsilon}{4}\cdot\frac{y}{\log y},$$
if $y$ is large enough.
Moreover,
\begin{multline*}
|\{pv\ |\ y^{2/3}<p\leq y/x,p \text{ is a prime, }v<\sqrt{x}\}|\leq \pi(y/x)\sqrt{x}\leq \\
\leq2\frac{y/x}{\log(y/x)}\sqrt{x}=2\frac{y}{\log y}\frac{1}{\sqrt{x}}\frac{1}{1-\frac{\log x}{\log y}}\leq \frac{\varepsilon}{4}\frac{y}{\log y},
\end{multline*}
since $x^4<y$ and $x\geq N>(32/(3\varepsilon))^2$.

Let us continue with the estimation of the next term:
\begin{multline*}
|\{pv\ |\ y/x<p\leq y/N,p\text{ is a prime}, v\leq\sqrt{ y/p}\}|\leq \\
\leq |\{pv\ |\ \exists j:\ N\leq j\leq x-1, y/(j+1)<p\leq y/j,p\text{ is a prime}, v\leq\sqrt{j+1}\}|\leq \\
\leq\sum\limits_{j=N}^{x-1} \left(\pi\left(\frac{y}{j}\right)-\pi\left(\frac{y}{j+1}\right)\right)\sqrt{j+1}
\end{multline*}
If $x$ is fixed and $y\to \infty$, then $\pi\left(\frac{y}{j}\right)=\frac{y}{j\log y}+o\left(\frac{y}{j(\log y)^{1.5}}\right)$,
therefore we obtain that
$$\pi\left(\frac{y}{j}\right)-\pi\left(\frac{y}{j+1}\right)=\frac{1}{j(j+1)}\frac{y}{\log y}+o\left(\frac{y}{j(\log y)^{1.5}}\right).$$
(For instance it suffices to take $y=\left\lfloor x^x \right\rfloor$.)
Hence,
$$\sum\limits_{j=N}^{x-1} \left(\pi\left(\frac{y}{j}\right)-\pi\left(\frac{y}{j+1}\right)\right)\sqrt{j+1}= \left(\sum\limits_{j=N}^{x-1}\frac{1}{j\sqrt{j+1}}\right)\frac{y}{\log y}+o\left(\frac{\sqrt{x}}{\sqrt{\log y}}\right)\frac{y}{\log y}.$$

Therefore,
$$|\{pv\ |\ y/x<p\leq y/N,p\text{ is a prime}, v\leq\sqrt{ y/p}\}|\leq \frac{\varepsilon}{4}\frac{y}{\log y}, $$
if $N>(16/\varepsilon)^2$ and $y$ is large enough.

Finally,
$$|\{p\ |\ y/N<p\leq y,p\text{ is a prime}\}|\leq \pi(y)\leq \left(1+ \frac{\varepsilon}{4}\right)\frac{y}{\log y},$$
if $y$ is large enough.
Adding up the estimates we obtain that
$$|B_{i+1}|\leq \left(1+\varepsilon\right)\frac{y}{\log y}$$
holds, if $y$ is sufficiently large. $\blacksquare $


The logarithmic density of a set $B\subset \mathbb{Z}^+$ is defined as the limit  $\displaystyle \lim_{n\to \infty }\frac{\sum\limits_{b\in B(n)}\frac{1}{b}}{\log n}$ (if it exists). Now, we prove Theorem~\ref{MBlog} which determines how small  $\sum\limits_{b\in B(n)}\frac{1}{b}$ can be for a multiplicative basis of order $h$.


\noindent
{\bf Proof of Theorem \ref{MBlog}}. In order to prove the first statement let $B\in MB_h(\mathbb{Z}^+)$, moreover let $B=\{b_1,b_2,\dots \}$, where $1\leq b_1<b_2<\dots$. Let us denote by $s_{B,h}(k)$ that how many ways $k$ can be written as a product of $h$ elements of the set $B$, that is,
$$s_{B,h}(k)=|\{ (i_1,\dots ,i_h)\in (\mathbb{Z}^+)^h:\ b_{i_1}\cdot\dots \cdot b_{i_h}=k\}|.$$
Clearly,
$$\left( \sum_{b\leq B(n)}\frac{1}{b}\right)^h\geq \sum_{k\leq n}\frac{s_{B,h}(k)}{k}\geq \sum_{k\leq n, |\mu(k)|=1}\frac{s_{B,h}(k)}{k},$$
where $\mu$ is the M\"obius-function, that is, the summation ranges over the squarefree integers.
If there exists a representation $k=b_{i_1}\dots b_{i_h}$, with $b_{i_1}<b_{i_2}<\dots <b_{i_h}$, then $s_{B,h}(k)\geq h!$ holds. On the other hand, if for some squarefree integer $k$ a representation $k=b_{i_1}\dots b_{i_h}$ with $b_{i_1}<b_{i_2}<\dots <b_{i_h}$ does not exist, then every representation of $k$ as a product of $h$ factors contains $b_1=1$. Hence,
$$\sum_{\substack{k\leq n,|\mu(k)|=1, \\ s_{B,h}(k)<h!}}\frac{s_{B,h}(k)}{k}\leq \sum_{\substack{k\leq n,|\mu(k)|=1, \\ s_{B,h}(k)<h!}}\frac{h!}{k} \leq  h!\left( \sum_{b\in B(n)}\frac{1}{b}\right) ^{h-1}.$$
Thus
\begin{multline*}
\left( \sum_{b\in B(n)}\frac{1}{b}\right)^h\geq \sum_{k\leq n, |\mu(k)|=1}\frac{s_{B,h}(k)}{k} = \\
=\sum_{\substack{k\leq n, |\mu(k)|=1, \\ s_{B,h}(k)\geq h!}}\frac{s_{B,h}(k)}{k}+\sum_{\substack{k\leq n, |\mu(k)|=1, \\ s_{B,h}(k)<h!}}\frac{s_{B,h}(k)}{k}\geq \\
\geq \sum_{\substack{k\leq n, |\mu(k)|=1, \\ s_{B,h}(k)\geq h!}}\frac{h!}{k}+\left( \sum_{\substack{k\leq n, |\mu(k)|=1, \\ s_{B,h}(k)< h!}}\frac{h!}{k} -h!\left( \sum_{b\in B(n)}\frac{1}{b}\right)^ {h-1} \right).
\end{multline*}
After some ordering we obtain by the binomial theorem that
$$\left( \left(\sum_{b\in B(n)}\frac{1}{b}\right) +(h-1)!\right) ^h \geq  \sum_{k\leq n, |\mu(k)|=1}\frac{h!}{k} .$$
Applying the well-known estimate $\sum\limits_{k\leq n, |\mu(k)|=1}\frac{1}{k}=\left( \frac{6}{\pi ^2}+o(1)\right) \log n $ and the inequalities $\sqrt{\frac{6}{\pi ^2}}\leq \sqrt[h]{\frac{6}{\pi ^2}}$ and $h!>\left( \frac{h}{e}\right )^h$ we obtain the bound claimed in the first part of the theorem.

To prove the second statement of the theorem let us denote the set of prime numbers by $P$ and the $k$th prime by $p_k$. First we show that the set $P$ can be partitioned into $h$ subsets, $P=P_1\cup \dots \cup P_h$, in such a way that
\begin{equation}\label{diff}
\max _{1\leq i\leq h}\left\{\ \sum_{p\in P_i(k)}\frac{1}{p}\right\}-\min _{1\leq i\leq h}\left\{\ \sum_{p\in P_i(k)}\frac{1}{p}\right\}\leq 0.5
\end{equation}
 hold for every $k\in \mathbb{Z}^+$. Let $P_i=\{ p_{i+hm}:m\ge 0\}$ for every $1\le i\le h$. Then it is easy to see that $\max\limits_{1\leq i\leq h}\left\{\ \sum\limits_{p\in P_i(k)}\frac{1}{p}\right\}=\sum\limits_{p\in P_1(k)}\frac{1}{p}$ and $\min\limits_{1\leq i\leq h}\left\{\ \sum\limits_{p\in P_i(k)}\frac{1}{p}\right\}=\sum\limits_{p\in P_h(k)}\frac{1}{p}$, moreover  $\sum\limits_{p\in P_h(k)}\frac{1}{p}>\sum\limits_{p\in P_1(k)\setminus \{ 2\}  }\frac{1}{p}=\left(\sum\limits_{p\in P_1(k)}\frac{1}{p}\right)-0.5$, which proves that the defined partition of the set of prime numbers satisfies \eqref{diff}.

For this partition $P=P_1\cup \dots \cup P_h$, we also have
$$\left|\sum _{p\in P_i( n)}\frac{1}{p}-\frac{1}{h}\sum _{p\in P(n)}\frac{1}{p}\right|\leq 0.5.$$
Now let us choose  the sets $A_i$ for $1\leq i\leq h$ in such  a way that in the set $A_i$ every integer $k$ has  prime factors only from the set $P_i$, that is,
$$A_i=\{k\in\mathbb{Z}^+:\text{each prime factor of $k$ belongs to $P_i$}\}.$$
Finally, let $\displaystyle B=\bigcup\limits_{i=1}^hA_i$. It is easy to see that $B$ is a multiplicative basis of order $h$. Therefore, it remains to prove that $\displaystyle \sum_{b\in B(n)}\frac{1}{b}\leq Ch\sqrt[h]{\log n}$ for some absolute constant $C$. Obviously,
$$\sum_{b\in B(n)}\frac{1}{b} =\sum_{i=1}^h\sum_{b\in A_i(n)}\frac{1}{b},$$
so it is enough to show that $\sum\limits_{b\in A_i(n)}\frac{1}{b}\leq C\sqrt[h]{\log n}$. Clearly,
$$\sum_{b\in A_i(n)}\frac{1}{b}\leq \prod _{p\in P_i(n)}\left( 1+\frac{1}{p}+\frac{1}{p^2}+\dots \right).$$
 Moreover, the following inequality holds for every prime number $p$:
$$1<\frac{1+\frac{1}{p}+\frac{1}{p^2}+\dots }{e^{\frac{1}{p}}}=\frac{1+\frac{1}{p-1}}{e^{\frac{1}{p}}}<e^{\frac{1}{p-1}-\frac{1}{p}}\leq e^{\frac{2}{p^2}}.$$
By the inequality $\sum\frac{1}{p^2}<\frac{1}{2} $ we obtain that
$$\prod _{p\in P_i(n)}\left( 1+\frac{1}{p}+\frac{1}{p^2}+\dots \right) \leq e\prod\limits_{p\in P_i(n)}e^{\frac{1}{p}} \leq e^{\frac{3}{2}}e^{\frac{1}{h}\sum\limits_{p\in P(n)}\frac{1}{p}},$$
 but the well-known estimate $\sum\limits_{p\in P(n)}\frac{1}{p}= \log \log n +O(1)$ gives that $e^{\frac{1}{h}\sum\limits_{p\in P(n)}\frac{1}{p}}=O(\sqrt[h]{\log n})$, which completes the proof. $\blacksquare $

\noindent
We continue with proving Theorem~\ref{rai2} which strengthens Raikov's result.

\noindent
{\bf Proof of Theorem \ref{rai2}}. To verify the lower bound $\displaystyle \limsup _{n\to \infty} \frac{|B(n)|}{n/\log ^{\frac{h-1}{h}} n }\ge \frac{\sqrt{6}}{e\pi }$, by Theorem~\ref{MBlog} it is enough to show  that  $\displaystyle \limsup _{n\to \infty}\frac{\sum\limits_{b\in B(n)}\frac{1}{b}}{h\sqrt[h]{\log n} }\le \limsup _{n\to \infty} \frac{|B(n)|}{n/\log ^{\frac{h-1}{h}} n }$. Let $B=\{ b_1,b_2,\dots \}$, where $0<b_1<b_2\dots $. If $\displaystyle \limsup _{n\to \infty} \frac{|B(n)|}{n/\log ^{\frac{h-1}{h}} n } =c$, then for every $\varepsilon >0$ there exists an $n_0=n_0(\varepsilon)$ such that $\displaystyle \frac{|B(n)|}{n/\log ^{\frac{h-1}{h}} n }<c+\varepsilon$ for $n\ge n_0$. Hence, $\displaystyle\frac{|B(b_n)|}{b_n/\log ^{\frac{h-1}{h}} b_n }<c+\varepsilon$, therefore $(c+\varepsilon) b_n>n\log ^{\frac{h-1}{h}} b_n \ge n\log ^{\frac{h-1}{h}} n$, that is, $\displaystyle b_n>\frac{1}{c+\varepsilon }n\log ^{\frac{h-1}{h}} n$ for $n\ge n_0$. Thus $\sum\limits_{b\in B(n)}\frac{1}{b} \le \sum\limits_{k=1}^n\frac{1}{b_k} = \sum\limits_{k<n_0}\frac{1}{b_k}+\sum\limits_{k=n_0}^n\frac{1}{b_k}<C_0+\sum\limits_{k=n_0}^n\frac{c+\varepsilon}{k\log ^{\frac{h-1}{h}} k}<C_1+(c+\varepsilon)\int\limits_{10}^n\frac{1}{x\log ^{\frac{h-1}{h}} x}dx<C_2+(c+\varepsilon )h\log ^{\frac{1}{h}} n$, which completes the proof of the first part.
\\
To prove the second statement it is enough to construct a multiplicative basis $B$ of order $h$ for which $\displaystyle \limsup _{n\to \infty} \frac{|B(n)|}{n/\log ^{\frac{h-1}{h}} n } =c_h$, where the sequence $c_h$ is bounded. We are going to show that Raikov's construction (see \cite{rai}) is a suitable choice for $B$. The set of prime numbers is denoted by $2=p_1<p_2<p_3<\dots $. The prime numbers are distributed into $h$ subsets in the following way: $P=\bigcup\limits_{i=1}^hP_i$, where $P_i=\{ p_{i+hm}:m\ge 0\}$. For $1\leq i\leq h$ let $N_i=\{a\in \mathbb{Z}^+: \text{each prime factor of $a$ belongs to $P_i$} \}$. We have already seen in the proof of Theorem~\ref{MBlog} that the set $N=\bigcup\limits_{i=1}^hN_i$ forms a multiplicative basis of order $h$. Let $\eta _i(s)=\prod\limits_{m=0}^{\infty}(1-p_{i+hm}^{-s})^{-1}$, if $s\in \mathbb{C}$ and $\Re(s)>1$. Let us fix the integer $1\le i\le h $ and take $\varphi _i(s)=\prod\limits_{q\ne i}\left( \frac{\eta _i(s)}{\eta _q(s)}\right) ^{1/h}$. Raikov proved that there exists an $\varepsilon >0$ such that the function $\varphi _i(s)$ is analytic for $\Re(s)>1-\varepsilon$ and $|N(x)|\sim \frac{\varphi _1(1)+\varphi _2(1)+\dots +\varphi _h(1)}{\Gamma \left(\frac{1}{k}\right)}x\log ^{\frac{1}{h}-1}x$. Since  $\displaystyle \Gamma \left( \frac{1}{h}\right) =\int _0^{\infty} x^{\frac{1}{h}-1}e^{-x}dx>\int _0^1 e^{-1}x^{\frac{1}{h}-1}dx=\frac{h}{e}$, therefore it is enough to prove that $\varphi _i(1)$ is bounded. Later on $s$ will denote a real number. We will show that for some suitable constants $0<c_2<1<c_3$ we have $\displaystyle c_2 <\lim _{s\to 1+}\frac{\eta _i(s)}{\eta _q(s)}<c_3$, therefore $c_2\le \varphi _i(1)\le c_3$. For $s>1$ we have that
\begin{multline*}
\frac{\eta _i(s)}{\eta _q(s)}=\prod _{m=0}^{\infty} \frac{1-p_{q+mh}^{-s}}{1-p_{i+mh}^{-s}}=\\
\exp\left(\sum_{m=0}^{\infty}\left( \log\left( 1-p_{q+mh}^{-s} \right) - \log\left( 1-p_{i+mh}^{-s} \right) \right)\right)=\exp
\left( -\sum_{m=0}^{\infty} \sum_{t=1}^{\infty} \left( \frac{p_{q+mh}^{-ts}}{t} - \frac{p_{i+mh}^{-ts}}{t} \right) \right).
\end{multline*}
A routine calculation gives that $\displaystyle  \left|\sum_{m=0}^{\infty} \sum_{t=2}^{\infty} \left( \frac{p_{q+mh}^{-ts}}{t} - \frac{p_{i+mh}^{-ts}}{t} \right)\right|<c_4$, therefore it remains to prove that $$c_5<\sum_{m=0}^{\infty}\left( \frac{1}{p_{q+mh}^{s}} - \frac{1}{p_{i+mh}^{s}} \right) <c_6.$$
Let us introduce a constant $x$ which will be defined later. Hence
\begin{multline*}
\left|\sum_{m=0}^{\infty}\left( \frac{1}{p_{q+mh}^{s}} - \frac{1}{p_{i+mh}^{s}} \right)\right|< \\
<\Big|\sum_{m:\ p_{q+mh}\le x}\left( \frac{1}{p_{q+mh}^{s}} - \frac{1}{p_{q+mh}} \right) -\sum_{m:\ p_{i+mh}\le x}\left( \frac{1}{p_{i+mh}^{s}} - \frac{1}{p_{i+mh}} \right) - \\
  \sum_{m:\ p_{q+mh}\le x,p_{i+mh}\le x }\left( \frac{1}{p_{i+mh}} - \frac{1}{p_{q+mh}} \right)\Big|+\sum _{p>x}\frac{1}{p^s}+c_7<\\
\sum_{m:\ p_{i+mh}\le x}\left( \frac{1}{p_{i+mh}} - \frac{1}{p_{i+mh}^{s}} \right) + \sum_{m:\ p_{q+mh}\le x}\left( \frac{1}{p_{q+mh}} - \frac{1}{p_{q+mh}^{s}} \right) \\
 +\left|\sum_{m:\ p_{i+mh}\le x, p_{q+mh}\le x}\left( \frac{1}{p_{i+mh}} - \frac{1}{p_{q+mh}} \right)\right| + \sum _{p>x}\frac{1}{p^s}+c_7.
\end{multline*}
The well-known estimation $1-y<e^{-y}$ yields $\frac{1}{p_{q+mh}} - \frac{1}{p_{q+mh}^{s}}<(s-1)\frac{\log p_{q+mh}}{p_{q+mh}}$, therefore $\sum\limits_{m:\ p_{q+mh}\le x}\left( \frac{1}{p_{q+mh}} - \frac{1}{p_{q+mh}^{s}}\right)\le (s-1)\sum\limits_{m:\ p_{q+mh}\le x} \frac{\log p_{q+mh}}{p_{q+mh}}\le c_8(s-1)\log x$. Similarly, $\sum\limits_{m:\ p_{i+mh}\le x}\left( \frac{1}{p_{i+mh}} - \frac{1}{p_{i+mh}^{s}}\right)\le c_8(s-1)\log x$.

We have seen in the proof of Theorem~\ref{MBlog} that $\left|\displaystyle \sum_{p\in P_u, p\le x}\frac{1}{p}- \frac{1}{h}\sum_{p\le x}\frac{1}{p}\right|\le 0.5$. Hence,
\begin{multline*}
\left|\sum_{m:\ p_{i+mh}\le x, p_{q+mh}\le x}\left( \frac{1}{p_{i+mh}} - \frac{1}{p_{q+mh}} \right)\right|= \\
\left|\left( \sum_{m:\ p_{i+mh}\le x, p_{q+mh}\le x}\frac{1}{p_{i+mh}}-\frac{1}{h}\sum_{p, p\le x}\frac{1}{p}\right)-\left( \sum_{m:\ p_{i+mh}\le x, p_{q+mh}\le x}\frac{1}{p_{q+mh}}-\frac{1}{h}\sum_{p, p\le x}\frac{1}{p}\right)\right|\le 1.
\end{multline*}

By the Prime Number Theorem we have
\begin{multline*}
\sum _{p>x}\frac{1}{p^s}\le \sum_{k\ge c_8\frac{x}{\log x}}\frac{1}{(k\log k)^s}\le \frac{c_9}{(\log x)^s}\sum_{k\ge c_8\frac{x}{\log x}}\frac{1}{k^s}\le \\
\frac{c_{10}}{(\log x)^s}\int\limits_{c_8\frac{x}{\log x}}^{\infty}\frac{dt}{t^s}\le \frac{c_{11}}{s-1}\frac{1}{x^{s-1}\log x}.
\end{multline*}
\noindent
Summarizing these bounds we get
$$\left|\sum_{m=0}^{\infty}\left( \frac{1}{p_{q+mh}^{s}} - \frac{1}{p_{i+mh}^{s}} \right)\right|<c_7+2c_{12}(s-1)\log x +\frac{c_{11}}{s-1}\frac{1}{x^{s-1}\log x}+1.$$
Substituting $x=e^{\frac{1}{s-1}}$ we get
$$\left|\sum_{m=0}^{\infty}\left( \frac{1}{p_{q+mh}^{s}} - \frac{1}{p_{i+mh}^{s}} \right)\right|<c_7+2c_{12}+\frac{c_{11}}{e}+1,$$ which completes the proof. $\blacksquare $

\noindent
The following lemma is going to be used in the proof of Theorem~\ref{Erdinf}:
\begin{lemma}\label{c+eps}
Let $Q$ be a subset of the prime numbers satisfying $|Q(n)|\ll n^c$ for some constant $c>0$. Then for every $\varepsilon>0 $ there exists some integer $N_0=N_0(\varepsilon ,Q)$ such that for every $n\geq N_0$ we have
$$|\{ k: k\leq n\hbox{ and every prime divisor of $k$ is in $Q$}\}|\leq n^{c+\varepsilon}.$$
\end{lemma}

\noindent
{\bf Proof of Lemma \ref{c+eps}}. Let the primes in $Q$ be: $q_1<q_2<\dots $ and denote by $p_n$ the $n$th prime number. Let us define an injective mapping $k\to k'$ in such a way that to $k=q_{i_1}^{\alpha _{i_1}}\dots q_{i_s}^{\alpha _{i_s}}$ we assign  $k'=p_{i_1}^{\alpha _{i_1}}\dots p_{i_s}^{\alpha _{i_s}}$. It is enough to prove that there exists a suitable set $Y$ satisfying $|Y|=(\log n)^{O(1)}$ such that each $k'$ can be represented as $k'=xy$, where $x\leq n^{c+\frac{\varepsilon}{2}}$ and $y\in Y$.  We know that $n=|Q(q_n)|\ll q_n^c$, 
hence $p_n\sim n\log n \ll q_n^c\log q_n$. Thus there exists some $C=C(\varepsilon, Q)$ such that for every $p_l>C$, we have $p_l\leq q_l^{c+\frac{\varepsilon}{2}}$. In $k'=\left(\prod\limits_{p_l\leq C}p_l^{\alpha _l}\right)\left(\prod\limits_{p_l> C}p_l^{\alpha _l}\right)$ we have $\prod\limits_{p_l> C}p_l^{\alpha _l}\leq n^{c+\frac{\varepsilon}{2}}$, if $k\leq n$. Furthermore,  for the product $\prod\limits_{p_l\leq C}p_l^{\alpha _l}$ there are at most  $(\log _2 n)^C$ possibilities, which completes the proof. $\blacksquare $


\smallskip

\noindent
Finally, we prove Theorem~\ref{Erdinf} about the infinite case of Erd\H os' problem.

\smallskip
\noindent
{\bf Proof of Theorem \ref{Erdinf}}. Let $P$ be the set of primes, moreover let $P_1=P\cap A$ and $P_2=P\setminus A$. Therefore, $P_1$ and $P_2$ are disjoint and $P_1\cup P_2=P$.
If $P_1=P$, then $A=P$, otherwise let $\displaystyle \alpha =\limsup _{n\to \infty}\frac{\log |P_2(n)|}{\log n}$. Let $B$ be a multiplicative basis of order $h$ defined in Lemma~\ref{basis}, and take the mapping $A(n)\to B$ defined in Lemma~\ref{inj}. We claim that if $a\to k$, where $k$ is not a prime number, then each prime factor of $k$ belongs to the set $P_2$.
Since, if $p|k$ for some $p\in P_1$, then for $a_0=p\in A$ and $a_1=a\in A$ we have $a_0|a_1$, which contradicts the assumption that $A$ possesses property $\mathcal{P}_h$.


If $\alpha =0$, then by Lemma~\ref{c+eps}:
\begin{multline*}
|A(n)|\leq \pi(n)+|\{ k: k\leq n^{\frac{2}{h+1}},\hbox{ $k$ is an image in the mapping $A\to B$}\} |\leq \\
\leq \pi (n)+|\{ k: k\leq n^{\frac{2}{h+1}}\hbox{ and each prime factor of $k$ belongs to the set $P_2$}\} |= \\
=\pi(n)+O\left( \left( n^{\frac{2}{h+1}}\right) ^{\varepsilon+\varepsilon/2} \right)  =\pi(n)+O\left( n^{\varepsilon }\right).
\end{multline*}

If $\alpha >0$, then we prove that $\displaystyle \liminf _{n\to \infty}(|A(n)|-\pi(n))=-\infty $. First we show that for every $ \delta >0 $ there exist infinitely many integers $n$ such that $|[2^{n-1},2^n]\cap P_2|>2^{(\alpha -\delta )n}$. For the sake of contradiction assume that $|[2^{n-1},2^n]\cap P_2|\leq 2^{(\alpha -\delta )n}$ for $n\geq n_0$, then $|P_2(2^N)|\leq M_0+1+2^{(\alpha -\delta )1}+2^{(\alpha -\delta )2}+\dots +2^{(\alpha -\delta )N}=O\left( 2^{(\alpha -\delta )N}\right) $. Hence for every $n$ we have $|P_2(n)|=O(n^{\alpha -\delta })$, which contradicts the definition of $\alpha$. By 
Lemma~\ref{basis} we have
\begin{multline*}
|A(2^n)|\leq |\{ p:p\leq 2^n\hbox{ and $p$ is an image in the mapping $A\to B$}\}|+ \\
+|\{ k: k\leq (2^n)^{\frac{2}{h+1}},\hbox{ $k$ is an image in the mapping $A\to B$}\} |\leq \\
\leq \pi (2^n)-|P_2\cap [2^{n-1},2^n]|+ |\{ k: k\leq (2^n)^{\frac{2}{h+1}},\hbox{ each prime factor of $k$ belongs to $P_2$}\} |
\end{multline*}
For infinitely many $n$ this can be bounded by
$$|A(2^n)|\leq \pi (2^n)-2^{(\alpha -\delta)n}+\left( (2^n)^{\frac{2}{h+1}} \right) ^{\alpha +\varepsilon}.$$
The values $\varepsilon =\delta =\frac{\alpha }{8}$ verify the desired statement.

In order to construct an always dense set $A$, a sequence $l_n\to\infty$ is going to be chosen and sequences $f_n$ and $g_n$ are going to be defined recursively as follows: Let $f_1$ and $g_1$ be large enough (we will specify them later), $f_{n+1}=\left( \frac{f_n}{2l_nh}\right) ^{l_n}$ and $g_{n+1}=g_n^{hl_n}$. Then it easy to see that
$$f_{n+1}=\left(  f_1\left( 2h \right) ^{-1-\frac{1}{l_1}-\frac{1}{l_1l_2}-\dots -\frac{1}{l_1l_2\dots l_{n-1}}} l_1^{-1}l_2^{-\frac{1}{l_1}}l_3^{-\frac{1}{l_1l_2}}\dots l_n^{-\frac{1}{l_1l_2\dots l_{n-1}}} \right) ^{l_1l_2\dots l_n}$$
and
 $$g_{n+1}=g_1^{h^nl_1l_2\dots l_n}.$$
 Let us suppose that for every $k$ we have $l_k\ge 2$ and $\displaystyle \sum _{n=1}^{\infty}\frac{\log l_n}{l_1l_2\dots l_{n-1}}<\infty$. In this case for some $c_1>1$ we have  $f_{n+1}>c_1^{l_1l_2\dots l_n}$, if $f_1$ is large enough. The set $A$ is defined with a little modification of the set of prime numbers. After the integer $g_m$ we just omit $f_m$ prime numbers $P_m=\{ p_1^{(m)},\dots ,p_{f_m}^{(m)}\}$, and instead of them we add the integers from the set
$$B_m=\left\{ p_{j_1^{(i)}}^{(m)}\dots p_{j_{hl_m-h+1}^{(i)}}^{(m)}: 1\leq i\leq \left(\frac{f_m}{2(hl_m-h+1)}\right)^{l_m} \right\},$$
where the sets $S_i=\{ j_1^{(i)},\dots ,j_{hl_m-h+1}^{(i)}  \}$ are the sets defined in Lemma~\ref{tinter} for $n=f_m$, $k=hl_m-h+1$ and $t=l_m$. If $g_1$ is large enough, then for every $m\geq 1$ the elements of $B_m$ are less than $g_{m+1}$. Moreover,  it is easy to check that $b_0\nmid b_1\dots b_h$ for any distinct integers from the set $B_m$,  therefore the set
$$A=\left( P\setminus  \bigcup\limits_{n=1}^{\infty}P_n\right) \cup \left( \bigcup\limits_{n=1}^{\infty} B_n \right)$$
possesses property $\mathcal{P}_h$. It remains to prove that the sequence $l_m$ can be chosen in such a way that
$$|A(x)|\geq \pi (x)+\exp\left\{{{(\log x)^{1-\frac{c\sqrt{\log h}}{\sqrt{\log \log x}}}}}\right\}$$
holds for every large enough $x$. Let us suppose that $g_{n+1}\leq x<g_{n+2}$. In this case
$$A(x)>\pi(x)-\left(\sum_{i=1}^{n+1}|P_i|\right)+|B_n|.$$
 An easy calculation gives that $\sum\limits_{i=1}^{n+1}|P_i|<1.1|P_{n+1}|$, if $n$ is large enough. On the other hand $l_n\to\infty$ implies that $|B_n|\ge \left( \frac{f_n}{2(hl_n-h+1)}\right) ^{l_n}\sim e^{\frac{h-1}{h}}\left( \frac{f_n}{2hl_n} \right)^{l_n}\sim e^{\frac{h-1}{h}}f_{n+1}=e^{\frac{h-1}{h}}|P_{n+1}|,$ therefore
$$A(x)-\pi (x)\gg f_{n+1}\gg c_1^{l_1\dots l_n}=\left( g_1^{h^{n+1}l_1\dots l_{n+1}}\right) ^{\frac{c_2}{h^{n+1}l_{n+1}}}=g_{n+2}^{\frac{c_2}{h^{n+1}l_{n+1}}}>x^{\frac{c_2}{h^{n+1}l_{n+1}}}.$$
It can be shown that an almost optimal (up to a constant factor) choice for $l_k$ is $l_k=h^k$. In this case $g_{n+1}\le x<g_{n+2}\textit{}$ can be rewritten as
$$g_1^{h^{n+\frac{n(n+1)}{2}}}\le x<g_1^{h^{n+1+\frac{(n+1)(n+2)}{2}}},$$
therefore $n\sim \sqrt{2\log _h \log x}$. Hence
 $$A(x)-\pi(x)>x^{c_2{h^{-\sqrt{c_3\log _h\log x}}}}>\exp\left\{{(\log x)^{1-\frac{c_4\sqrt{\log h}}{\sqrt{\log \log x}}}}\right\}. \blacksquare $$



\section{Questions}
Finally, we present some open problems.

\begin{problem}
Do there exist constants $c_h$ such that for the size of a smallest multiplicative basis for $[n]$ we have
$$ \displaystyle \min_{B\in MB_h([n])}|B|= \pi(n)+(c_h+o(1))\frac{n^{2/(h+1)}}{\log ^2 n}?$$
If so, determine $c_h$.
\end{problem}
\noindent
A similar problem can be formulated for $F_h(n)$.
\begin{problem}
Do there exist constants $d_h$ such that $$F_h(n)= \pi(n)+(d_h+o(1))\frac{n^{2/(h+1)}}{\log^2 n}?$$
\end{problem}
\noindent
We can not improve the lower bound in Theorem~\ref{Erdinf}. Is it true that it is almost optimal, that is:
\begin{problem}
Is it true that there exists constant $C_h$ such that if an infinite set $A$ satisfies $|A(n)|\geq \pi (n)+\exp\left\{{(\log n)^{1-\frac{C_h\sqrt{\log h}}{\sqrt{\log \log n}}}}\right\}$ for every $n$, then there are distinct elements $a_0,a_1,\dots, a_h\in A$ with $a_0| a_1\dots a_h$?
\end{problem}
\noindent
Let us denote by $F_{r,s}(n)$ the maximal size of $A\subset [n]$ such that there are no distinct elements $a_1,a_2,\dots ,a_{r+s}\in A$ with $a_1\dots a_r\mid a_{r+1}\dots a_{r+s}$. Clearly $F_{1,s}(n)=F_s(n)$. It is easy to see that for $1\leq r<s$ we have $F_{1,s}(n)\le F_{r,s}(n)\leq F_{s,s}(n)$. We know that $F_{1,3}(n)=\pi (n)+n^{1/2+o(1)}$ and following the estimation of the size of a multiplicative 3-Sidon sequence in \cite{ppp} it can be shown that $F_{3,3}(n)\leq \pi(n)+\pi(n/3)+n^{2/3+o(1)}$. Moreover, $F_{2,3}(n)\leq \pi(n)+n^{2/3+o(1)}$ can be deduced also.
\begin{problem}
Is it true that there exists a constant $c_{2,3}$ such that $F_{2,3}(n)=\pi (n)+n^{c_{2,3}+o(1)}$? If so, determine $c_{2,3}$.
\end{problem}

\end{document}